\documentclass[12pt]{article}
\usepackage{amsxtra}
\usepackage{amssymb, amsmath}
\usepackage{amsthm}
\usepackage {hyperref}

\newtheorem{theorem}{Theorem}[section]
\newtheorem{remark}{Remark}[section]

\newtheorem{proposition}{Proposition}[section]
\newtheorem{example}{Example}[section]

\begin{document}
\begin{center}
\textbf{\LARGE{Bounds on Non-Symmetric Divergence Measures in Terms
of Symmetric Divergence Measures}}
\end{center}

\begin{center}
\textbf{Inder Jeet Taneja}
\end{center}
\begin{center}
Departamento de Matem\'{a}tica\\
Universidade Federal de Santa Catarina\\
88.040-900 Florian\'{o}polis, SC, Brazil.\\
\textit{e-mail: taneja@mtm.ufsc.br\\
http://www.mtm.ufsc.br/$\sim $taneja}
\end{center}

\begin{abstract}
There are many information and divergence measures exist in the
literature on information theory and statistics. The most famous
among them are Kullback-Leibler \cite{kul} \textit{relative
information} and Jeffreys \cite{jef} \textit{J- divergence}. Sibson
\cite{sib} \textit{Jensen-Shannon divergence} has also found its
applications in the literature. The author \cite{tan3} studied a new
divergence measures based on \textit{arithmetic and geometric
means}. The measures like \textit{harmonic mean divergence} and
\textit{triangular discrimination} \cite{dac} are also known in the
literature. Recently, Dragomir et al. \cite{dgp} also studies a new
measure similar to \textit{J-divergence}, we call here the
\textit{relative J-divergence}. Another measures arising due to
Jensen-Shannon divergence is also studied by Lin \cite{lin}. Here we
call it \textit{relative Jensen-Shannon divergence}.
\textit{Relative arithmetic-geometric divergence} (Taneja
\cite{tan8}) is also studied here. All these measures can be written
as particular cases of Csisz\'{a}r \textit{f-divergence}. By putting
some conditions on the probability distribution, the aim here is to
obtain bounds among the measures.
\end{abstract}


\section{Introduction}

Let
\[
\Gamma _n = \left\{ {P = (p_1 ,p_2 ,...,p_n )\left| {p_i >
0,\sum\limits_{i = 1}^n {p_i = 1} } \right.} \right\}, \quad n
\geqslant 2,
\]

\noindent be the set of all complete finite discrete probability
distributions. There are many information and divergence measures
exists in the literature on information theory and statistics.
Some of them are symmetric with respect to probability
distributions, while others are not. Here we have divided these
measure in these two categories. Through out the paper it is under
stood that the probability distributions $P,Q \in \Gamma _n $.

\subsection{Non-Symmetric Measures }
Here we shall give some non-symmetric measures of information. The
most famous among them are $\chi ^2 - $divergence and
Kullback-Leibler relative information. We understand by non
symmetric measures are those that are not symmetric with respect to
probability distributions $P,Q \in \Gamma _n $. These measures as
follows.\\

\textbf{$\bullet$ $\chi ^2 - $Divergence} (Pearson \cite{pea})
\begin{equation}
\label{eq1}
\chi ^2(P\vert \vert Q) = \sum\limits_{i = 1}^n {\frac{(p_i - q_i )^2}{q_i
}} = \sum\limits_{i = 1}^n {\frac{p_i^2 }{q_i } - 1}
\end{equation}

\textbf{$\bullet$ Relative Information} (Kullback and Leibler
\cite{kul})
\begin{equation}
\label{eq3}
K(P\vert \vert Q) = \sum\limits_{i = 1}^n {p_i \ln (\frac{p_i }{q_i })}
\end{equation}

\textbf{$\bullet$ Relative J-Divergence} (Dragomir et al.
\cite{dgp})
\begin{equation}
\label{eq5}
D(P\vert \vert Q) = \sum\limits_{i = 1}^n {(p_i - q_i )\ln \left( {\frac{p_i
+ q_i }{2q_i }} \right)}
\end{equation}

\textbf{$\bullet$ Relative Jensen-Shannon divergence} (Sibson
\cite{sib})
\begin{equation}
\label{eq7}
F(P\vert \vert Q) = \sum\limits_{i = 1}^n {p_i \ln \left( {\frac{2p_i }{p_i
+ q_i }} \right)}
\end{equation}

\textbf{$\bullet$ Relative arithmetic-geometric divergence} (Taneja
\cite{tan3})
\begin{equation}
\label{eq9}
G(P\vert \vert Q) = \sum\limits_{i = 1}^n {\left( {\frac{p_i + q_i }{2}}
\right)\ln \left( {\frac{p_i + q_i }{2p_i }} \right)}
\end{equation}

\subsection{Symmetric Measures of Information }
Here we shall give some symmetric measures of information. Some of
them can be obtained from subsection 1.1. These measures as
follows.\\

\textbf{$\bullet$ Hellinger Discrimination} (Hellinger \cite{hel})
\begin{equation}
\label{eq11} h(P\vert \vert Q) = 1 - B(P\vert \vert Q) =
\frac{1}{2}\sum\limits_{i = 1}^n {(\sqrt p_i - \sqrt q_i )^2} .
\end{equation}

\noindent where
\begin{equation}
\label{eq12} B(P\vert \vert Q) = \sum\limits_{i = 1}^n{ \sqrt {p_i
q_i }}
\end{equation}

\noindent is the well-known Bhattacharyya \cite{bha}
\textit{distance}.\\

\textbf{$\bullet$ Triangular Discrimination} (Dacunha-Castelle
\cite{dac})
\begin{equation}
\label{eq13}
\Delta (P\vert \vert Q) = 2\left[ {1 - W(P\vert \vert Q)} \right] =
\sum\limits_{i = 1}^n {\frac{(p_i - q_i )^2}{p_i + q_i }} .
\end{equation}

\noindent where
\begin{equation}
\label{eq14}
W(P\vert \vert Q) = \sum\limits_{i = 1}^n {\frac{2p_i q_i }{p_i + q_i }} .
\end{equation}

\noindent is the well-known \textit{harmonic mean divergence}.\\

\textbf{$\bullet$ Symmetric Chi-square Divergence} (Dragomir et al.
\cite{dsb})
\begin{align}
\Psi (P\vert \vert Q) & = \chi ^2(P\vert \vert Q) + \chi ^2(Q\vert
\vert P) \notag\\
\label{eq15} & = \sum\limits_{i = 1}^n {\frac{(p_i - q_i )^2(p_i +
q_i )}{p_i q_i }} .
\end{align}

\textbf{$\bullet$ J-divergence} (Jeffreys \cite{jef}; Kullback and
Leibler \cite{kul})
\begin{align}
 J(P\vert \vert Q) & = K(P\vert \vert Q) + K(Q\vert
\vert P)\notag\\
& = D(P\vert \vert Q) + D(Q\vert \vert P)\notag\\
\label{eq16} & = \sum\limits_{i = 1}^n {(p_i - q_i )\ln (\frac{p_i
}{q_i })}.
\end{align}

\textbf{$\bullet$ Jensen-Shannon divergence} (Sibson \cite{sib};
Burbea and Rao \cite{bra1})
\begin{align}
I(P\vert \vert Q) & = \frac{1}{2}\left[ {F(P\vert \vert Q) +
F(Q\vert \vert
P)} \right].\notag\\
&\label{eq17} = \frac{1}{2}\left[ {\sum\limits_{i = 1}^n {p_i \ln
\left( {\frac{2p_i }{p_i + q_i }} \right) + } \sum\limits_{i = 1}^n
{q_i \ln \left( {\frac{2q_i }{p_i + q_i }} \right)} } \right].
\end{align}

\textbf{$\bullet$ Arithmetic-Geometric divergence} (Taneja
\cite{tan3})
\begin{align}
T(P\vert \vert Q) & = \frac{1}{2}\left[ {G(P\vert \vert
Q) + G(Q\vert \vert P)} \right].\notag\\
& \label{eq18} = \sum\limits_{i = 1}^n {\left( {\frac{p_i + q_i
}{2}} \right)\ln \left( {\frac{p_i + q_i }{2\sqrt {p_i q_i } }}
\right)}.
\end{align}

After simplification, we can write
\[
J(P\vert \vert Q) = 4\left[ {I(P\vert \vert Q) + T(P\vert \vert Q)}
\right]
\]

\noindent and
\[
D(Q\vert \vert P) = \frac{1}{2}\left[ {F(P\vert \vert Q) + G(P\vert
\vert Q)} \right].
\]

\bigskip
The measures $I(P\vert \vert Q)$, $J(P\vert \vert Q), T(P\vert \vert
Q)$, $D(P\vert \vert Q), F(P\vert \vert Q)$ and $G(P\vert \vert Q)$
can be written in terms of $K(P\vert \vert Q)$ as follows:
\begin{align}
I(P\vert \vert Q) & = \frac{1}{2}\left[ {K\left( {P\vert \vert
\frac{P + Q}{2}} \right) + K\left( {Q\vert \vert \frac{P +
Q}{2}} \right)} \right],\notag\\
J(P\vert \vert Q) & = K(P\vert \vert Q) + K(Q\vert
\vert P),\notag\\
T(P\vert \vert Q) & = \frac{1}{2}\left[ {K\left( {\frac{P +
Q}{2}\vert \vert P} \right) + K\left( {\frac{P +
Q}{2}\vert \vert Q} \right)} \right],\notag\\
D(P\vert \vert Q) & = \frac{1}{2}\left[ {K\left( {Q\vert \vert
\frac{P + Q}{2}} \right) + K\left( {\frac{P + Q}{2}\vert \vert Q}
\right)} \right]\notag,\\
F(P\vert \vert Q) & = K\left( {P\vert \vert \frac{P + Q}{2}}
\right)\notag \intertext{and} G(P\vert \vert Q) & = K\left( {\frac{P
+ Q}{2}\vert \vert P} \right).\notag
\end{align}

\noindent respectively.

The following \textit{parallelogram identity} is also famous in the
literature \cite{csk}:
\[
K(P\vert \vert U) + K(Q\vert \vert U) = K\left( {P\vert \vert
\frac{P + Q}{2}} \right) + K\left( {Q\vert \vert \frac{P + Q}{2}}
\right) + 2\, K\left( {\frac{P + Q}{2}\vert \vert U} \right),
\]

\noindent for all $P,Q,U \in \Gamma _n $

Some studies on information and divergence measures can be seen in
Taneja \cite{tan1}, \cite{tan2}, \cite{tan3}. Also see on line
book by Taneja \cite{tan4}.

From the symmetric measures we observe that the measure
(\ref{eq12}) is a part of measure (\ref{eq11}) and the measure
(\ref{eq14}) is a part of measure (\ref{eq13}). Thus we have the
six measures (\ref{eq11}), (\ref{eq13}), (\ref{eq15})-(\ref{eq18})
symmetric with respect to probability distributions.

The following inequalities are already known:
\begin{equation}
\label{eq270} \frac{1}{2}h(P\vert \vert Q) \leqslant
\frac{1}{4}\Delta (P\vert \vert Q) \leqslant h(P\vert \vert Q),
\end{equation}
\begin{equation}
\label{eq271} \frac{1}{4}\Delta (P\vert \vert Q) \leqslant
h(P\vert \vert Q) \leqslant \frac{1}{16}\Psi (P\vert \vert Q)
\end{equation}
\begin{equation}
\label{eq272} \Delta (P\vert \vert Q) \leqslant \frac{1}{2}J(P\vert
\vert Q) \leqslant \frac{1}{4}\Psi (P\vert \vert Q)
\end{equation}

\noindent and
\begin{equation}
\label{eq273} \frac{1}{4}\Delta (P\vert \vert Q) \leqslant
I(P\vert \vert Q) \leqslant \frac{\log 2}{2}\Delta (P\vert \vert
Q).
\end{equation}

The inequalities (\ref{eq270}) are due to LeCam \cite{lec} and
Dacunha-Castelle \cite{dac}. The inequalities (\ref{eq271}) are due
to Taneja \cite{tan7}. The inequalities (\ref{eq272}) are due to
Dragomir et al. \cite{dsb}. Finally, the inequalities (\ref{eq273})
are due to Tops$\o$e \cite{top}.

Recently, Taneja \cite{tan9, tan10} proved the following
inequalities:
\begin{equation}
\label{eq32}
\frac{1}{4}\Delta (P\vert \vert Q) \leqslant I(P\vert \vert Q) \leqslant
h(P\vert \vert Q) \leqslant \frac{1}{8}J(P\vert \vert Q) \leqslant T(P\vert
\vert Q) \leqslant \frac{1}{16}\Psi (P\vert \vert Q).
\end{equation}

In this paper, our aim is to relate the \textit{non-symmetric
divergence measures} with the \textit{symmetric measures} given by
(\ref{eq13}), (\ref{eq16})-({\ref{eq18}). In order to obtain these
relationship we shall use the idea of \textit{Csisz\'{a}r
f-divergence} and in some cases making restrictions on the
probability distributions.

\section{$f-$Divergence and Information Measures}

Given a convex function $f:(0,\infty ) \to \mathbb{R}$, the $f -
$divergence measure introduced by Csisz\'{a}r \cite{csi1} is given
by
\begin{equation}
\label{eq33}
C_f (P\vert \vert Q) =
\sum\limits_{i = 1}^n {q_i f\left( {\frac{p_i }{q_i }} \right)} ,
\end{equation}

\noindent where $P,Q \in \Gamma _n $.

The following theorem is well known in the literature (ref.
Csisz\'{a}r \cite{csi1, csi2}).

\begin{theorem} \label{the21} Let the function $f:(0,\infty ) \to
\mathbb{R}$ is differentiable convex and normalized, i.e., $f(1) =
0$, then the Csisz\'{a}r $f - $divergence, $C_f (P\vert \vert Q)$,
given by (\ref{eq33}) is nonnegative and convex in the pair of
probability distribution $(P,Q) \in \Gamma _n \times \Gamma _n $.
\end{theorem}

The following theorem is due to Dragomir \cite{dra1, dra2}. It gives
bounds on \textit{Csisz\'{a}r f - divergence}.

\begin{theorem} \label{the22} (Dragomir \cite{dra1,
dra2}). Let $f:\mathbb{R}_ + \to \mathbb{R}$ be differentiable
convex and normalized i.e., $f(1) = 0$. Then
\begin{equation}
\label{eq24} 0 \leqslant C_f (P\vert \vert Q) \leqslant E_{C_f }
(P\vert \vert Q)
\end{equation}

\noindent where
\begin{equation}
\label{eq25} E_{C_f } (P\vert \vert Q) = \sum\limits_{i = 1}^n {(p_i
- q_i )} {f}'(\frac{p_i }{q_i }),
\end{equation}

\noindent for all $P,Q \in \Gamma _n $.

Let $P,Q \in \Gamma _n $ be such that there exists $r,R$ with $0 < r
\leqslant \frac{p_i }{q_i } \leqslant R < \infty $, $\forall i \in
\{1,2,...,n\}$, then
\begin{equation}
\label{eq26} 0 \leqslant C_f (P\vert \vert Q) \leqslant A_{C_f }
(r,R),
\end{equation}

\noindent where
\begin{equation}
\label{eq27} A_{C_f } (r,R) = \frac{1}{4}(R - r)\left[ {{f}'(R) -
{f}'(r)} \right].
\end{equation}

Further, if we suppose that $0 < r \leqslant 1 \leqslant R < \infty
$, $r \ne R$, then
\begin{equation}
\label{eq28} 0 \leqslant C_f (P\vert \vert Q) \leqslant B_{C_f }
(r,R),
\end{equation}

\noindent where
\begin{equation}
\label{eq29} B_{C_f } (r,R) = \frac{(R - 1)f(r) + (1 - r)f(R)}{R -
r}.
\end{equation}
\end{theorem}

Moreover, the following inequalities hold:
\begin{equation}
\label{eq30} E_{C_f } (P\vert \vert Q) \leqslant A_{C_f } (r,R),
\end{equation}
\begin{equation}
\label{eq31} B_{C_f } (r,R) \leqslant A_{C_f } (r,R)
\end{equation}

\noindent and
\begin{equation}
\label{eq32} 0 \leqslant B_{C_f } (r,R) - C_f (P\vert \vert Q)
\leqslant A_{C_f } (r,R).
\end{equation}

The inequalities (\ref{eq30}) and (\ref{eq32}) can be seen in
Dragomir \cite{dra2}, while the inequality (\ref{eq31}) can be
proved easily.

The following theorem is due to Taneja \cite{tan9, tan10}. It
relates two $f-$divergence measures.

\begin{theorem} \label{the23} Let $f_1 ,f_2 :I \subset \mathbb{R}_ + \to
\mathbb{R}$ be two differentiable convex functions which are
normalized, i.e., $f_1 (1) = f_2 (1) = 0$ and suppose that:

(i) $f_1 $and $f_2 $ are twice differentiable on $(r,R)$;

(ii) there exists the real constants $m, M$ such that $m < M$ and
\[
m \leqslant \frac{f_1 ^{\prime \prime }(x)}{f_2 ^{\prime \prime
}(x)} \leqslant M, \,\, f_2 ^{\prime \prime }(x)
> 0, \,\, \forall x \in (r,R)
\]
\noindent then we have
\begin{equation}
\label{eq43} m\,C_{f_2 } (P\vert \vert Q) \leqslant C_{f_1 } (P\vert
\vert Q) \leqslant M\,C_{f_2 } (P\vert \vert Q).
\end{equation}
\end{theorem}

\section{Bounds on Divergence Measures}

Based on Theorems \ref{the21} and \ref{the22}, we have the
particular cases for the measures given in Section 1. These
particular cases are given as examples, where the following the
expression is frequently used:
\begin{equation}
\label{eq47} L_{ - 1}^{ - 1} (a,b) = \begin{cases}
 {\frac{\ln b - \ln a}{b - a},} &{a \ne b} \\
 a & {a = b} \\
\end{cases}
\end{equation}
for all $a > 0$, $b > 0$.

\begin{example} \label{ex31} (\textit{Relative J-Divergence}). Let
us consider
\begin{equation}
\label{eq48} f_{D } (x) = (x - 1)\ln \left( {\frac{x + 1}{2}}
\right), \quad x \in (0,\infty )
\end{equation}

\noindent in (\ref{eq33}), then we have $C_f (P\vert \vert Q) =
D(P\vert \vert Q)$.

Moreover,
\begin{equation}
\label{eq49} {f}'_{D } (x) = \frac{x - 1}{x + 1} + \ln \left(
{\frac{x + 1}{2}} \right)
\end{equation}

\noindent and
\begin{equation}
\label{eq50} {f}''_{D } (x) = \frac{x + 3}{(x + 1)^2}.
\end{equation}

In view of (\ref{eq48}), (\ref{eq49}), Theorems \ref{the21} and
\ref{the22}, we have the following bounds on \textit{relative
J-divergence}:
\begin{equation}
\label{eq51} 0 \leqslant D(P\vert \vert Q) \leqslant E_{D } (P\vert
\vert Q) \leqslant A_{D } (r,R)
\end{equation}

\noindent and
\begin{equation}
\label{eq52} 0 \leqslant D(P\vert \vert Q) \leqslant B_{D } (r,R)
\leqslant A_{D } (r,R),
\end{equation}

\noindent where
\begin{align}
E_{D } (P\vert \vert Q) & = D(P\vert \vert Q) + \Delta (P\vert \vert
Q), \notag\\
A_{D } (r,R) & = \frac{1}{4}(R - r)^2\left[ {\frac{2}{(R + 1)(r +
1)} + L_{ - 1}^{ - 1} (r + 1,R + 1)} \right]\notag\\
\intertext{and} B_{D } (r,R) & = (R - 1)(1 - r)L_{ - 1}^{ - 1} (r +
1,R + 1).\notag
\end{align}
\end{example}

\begin{example} \label{ex33}(\textit{Relative Jensen-Shannon divergence}).
Let us consider
\begin{equation}
\label{eq64} f_{F } (x) = x\,\ln \left( {\frac{2x}{x + 1}} \right) -
\frac{x - 1}{2}, \quad x \in (0,\infty )
\end{equation}

\noindent in (\ref{eq33}), then we have $C_f (P\vert \vert Q) =
F(P\vert \vert Q)$.

Moreover,
\begin{equation}
\label{eq65} {f}'_{F } (x) = \frac{1}{2}\frac{x - 1}{x + 1} + \ln
\left( {\frac{2x}{x + 1}} \right)
\end{equation}

\noindent and
\begin{equation}
\label{eq66}
{f}''_{F_2 } (x) = \frac{1}{x(x + 1)^2}.
\end{equation}

In view of (\ref{eq64}), (\ref{eq65}), Theorems \ref{the21} and
\ref{the22}, we have the following bounds on \textit{relative
Jensen-Shannon divergence}:
\begin{equation}
\label{eq67} 0 \leqslant F(P\vert \vert Q) \leqslant E_{F } (P\vert
\vert Q) \leqslant A_{F } (r,R)
\end{equation}

\noindent and
\begin{equation}
\label{eq68} 0 \leqslant F(P\vert \vert Q) \leqslant B_{F } (r,R)
\leqslant A_{F } (r,R),
\end{equation}

\noindent where
\begin{align}
 E_{F } (P\vert \vert Q) & = D(Q\vert \vert P) -
\frac{1}{2}\Delta (P\vert \vert Q),\notag\\
A_{F } (r,R) & = \frac{1}{4}\frac{(R - r)^2}{(R + 1)(r + 1)}\left[
{L_{ - 1}^{ - 1} \left( {\frac{r}{r + 1},\frac{R}{R + 1}} \right) -
1} \right]\notag\\
\intertext{and} B_{F_2 } (r,R) & = \frac{1}{(R - r)}\left[ {R\ln
\left( {\frac{2R}{R + 1}} \right) - r\ln \left( {\frac{2r}{r + 1}}
\right)} \right]\notag\\
& \qquad \qquad - \frac{rR}{(R + 1)(r + 1)}L_{ - 1}^{ - 1} \left(
{\frac{r}{r + 1},\frac{R}{R + 1}} \right).\notag
\end{align}
\end{example}

\begin{example} \label{ex35}(\textit{Relative arithmetic-geometric divergence}). Let us consider
\begin{equation}
\label{eq80} f_{G } (x) = \frac{x + 1}{2}\,\ln \left( {\frac{x +
1}{2x}} \right) + \frac{x - 1}{2}, \quad x \in (0,\infty )
\end{equation}

\noindent in (\ref{eq33}), then we have $C_f (P\vert \vert Q) =
G(P\vert \vert Q)$.

Moreover,
\begin{equation}
\label{eq81} {f}'_{G } (x) = \frac{1}{2}\left[ {\ln \left( {\frac{x
+ 1}{2x}} \right) - \frac{x - 1}{x}} \right]
\end{equation}

\noindent and
\begin{equation}
\label{eq82}
{f}''_{G_2 } (x) = \frac{1}{2x^2(x + 1)}.
\end{equation}

In view of (\ref{eq80}), (\ref{eq81}), Theorems \ref{the21} and
\ref{the22}, we have the following bounds on \textit{relative
arithmetic-geometric divergence}:
\begin{equation}
\label{eq83} 0 \leqslant G(P\vert \vert Q) \leqslant E_{G } (P\vert
\vert Q) \leqslant A_{G } (r,R)
\end{equation}

\noindent and
\begin{equation}
\label{eq84} 0 \leqslant G(P\vert \vert Q) \leqslant B_{G } (r,R)
\leqslant A_{G } (r,R),
\end{equation}

\noindent where

\begin{align}E_{G } (P\vert \vert Q) & =
\frac{1}{2}\left[ {\chi ^2(Q\vert \vert P) - D(Q\vert \vert P)}
\right].\notag\\
A_{G } (r,R) & = \frac{(R - r)^2}{8rR}\left[ {1 - L_{ - 1}^{ - 1}
\left( {\frac{r + 1}{r},\frac{R + 1}{R}} \right)} \right]\notag\\
\intertext{and} B_{G } (r,R) & = \frac{1}{2}\ln \left[ {\frac{(r +
1)(R + 1)}{4rR}} \right] - \frac{1 - Rr}{2rR}L_{ - 1}^{ - 1} \left(
{\frac{r + 1}{r},\frac{R + 1}{R}} \right)\notag
\end{align}
\end{example}

\begin{example} \label{ex37} (\textit{Triangular discrimination}). Let us
consider
\begin{equation}
\label{eq96}
f_\Delta (x) = \frac{(x - 1)^2}{x + 1},
\quad
x \in (0,\infty )
\end{equation}

\noindent in (\ref{eq33}), then we have $C_f (P\vert \vert Q) =
\Delta (P\vert \vert Q)$.

Moreover,
\begin{equation}
\label{eq97}
{f}'_\Delta (x) = \frac{(x - 1)(x + 3)}{(x + 1)^2}
\end{equation}

\noindent and
\begin{equation}
\label{eq98}
{f}''_\Delta (x) = \frac{8}{(x + 1)^3}.
\end{equation}

In view of (\ref{eq96}), (\ref{eq97}), Theorems \ref{the21} and
\ref{the22}, we have the following bounds on \textit{triangular
discrimination}:
\begin{equation}
\label{eq99} 0 \leqslant \Delta (P\vert \vert Q) \leqslant E_{\Delta
} (P\vert \vert Q) \leqslant A_{\Delta } (r,R)
\end{equation}

\noindent and
\begin{equation}
\label{eq100} 0 \leqslant \Delta(P\vert \vert Q) \leqslant B_{\Delta
} (r,R) \leqslant A_{\Delta } (r,R),
\end{equation}

\noindent where
\begin{align}
E_\Delta (P\vert \vert Q) & = \sum\limits_{i = 1}^n {\left(
{\frac{p_i - q_i }{p_i + q_i }} \right)^2(p_i + 3q_i )} ,\notag\\
A_\Delta (r,R) & = \frac{(R - r)^2(R + r + 2)}{(R + 1)^2(r +
1)^2}\notag\\
\intertext{and} B_\Delta (r,R) & = \frac{2(R - 1)(1 - r)}{(R + 1)(1
+ r)}.\notag
\end{align}
\end{example}

\begin{example} \label{ex38} (\textit{J-divergence}). Let us consider
\begin{equation}
\label{eq104}
f_J (x) = (x - 1)\ln x,
\quad
x \in (0,\infty ),
\end{equation}

\noindent in (\ref{eq33}), then we have $C_f (P\vert \vert Q) =
J(P\vert \vert Q).$

Moreover,
\begin{equation}
\label{eq105}
{f}'_J (x) = 1 - x^{ - 1} + \ln x,
\end{equation}

\noindent and
\begin{equation}
\label{eq106}
{f}''_J (x) = \frac{x + 1}{x^2}.
\end{equation}

In view of (\ref{eq104}), (\ref{eq105}), Theorems \ref{the21} and
\ref{the22}, we have the following bounds on \textit{J-divergence}:
\begin{equation}
\label{eq107} 0 \leqslant J(P\vert \vert Q) \leqslant E_{J } (P\vert
\vert Q) \leqslant A_{J } (r,R)
\end{equation}

\noindent and
\begin{equation}
\label{eq108} 0 \leqslant J(P\vert \vert Q) \leqslant B_{J } (r,R)
\leqslant A_{J } (r,R),
\end{equation}

\noindent where
\begin{align}
E_J (P\vert \vert Q) & = J(P\vert \vert Q) + \chi ^2(Q\vert \vert
P),\notag\\
A_J (r,R) & = \frac{1}{4}(R - r)^2\left[ {(rR)^{ - 1} + L_{ - 1}^{ -
1} (r,R)} \right]\notag\\
\intertext{and} B_J (r,R) & = (R - 1)(1 - r)L_{ - 1}^{ - 1}
(r,R).\notag
\end{align}
\end{example}

\begin{example} \label{ex39} (\textit{Jensen-Shannon divergence}). Let us
consider
\begin{equation}
\label{eq112}
f_I (x) = \frac{x}{2}\ln x + \frac{x + 1}{2}\ln \left( {\frac{2}{x + 1}}
\right),
\quad
x \in (0,\infty ),
\end{equation}

\noindent in (\ref{eq33}), then we have $C_f (P\vert \vert Q) =
I(P\vert \vert Q).$

Moreover,
\begin{equation}
\label{eq113}
{f}'_I (x) = \frac{1}{2}\ln \left( {\frac{2x}{x + 1}} \right),
\end{equation}

\noindent and
\begin{equation}
\label{eq114}
{f}''_I (x) = \frac{1}{2x(x + 1)}.
\end{equation}

In view of (\ref{eq112}), (\ref{eq113}), Theorems \ref{the21} and
\ref{the22}, we have the following bounds on \textit{Jensen-Shannon
divergence}:
\begin{equation}
\label{eq115} 0 \leqslant I(P\vert \vert Q) \leqslant E_{I } (P\vert
\vert Q) \leqslant A_{I } (r,R)
\end{equation}

\noindent and
\begin{equation}
\label{eq116} 0 \leqslant I(P\vert \vert Q) \leqslant B_{I } (r,R)
\leqslant A_{I } (r,R),
\end{equation}
\noindent where
\begin{align}
E_I (P\vert \vert Q) & = \frac{1}{2}D(Q\vert \vert P), \notag\\
A_I (r,R) & = \frac{1}{8}\frac{(R - r)^2}{(R + 1)(r + 1)}L_{ - 1}^{
- 1} \left( {\frac{r}{r + 1},\frac{R}{R + 1}} \right)\notag\\
\intertext{and} B_I (r,R) & = \frac{1}{2(R - r)}\left[ {(R -
1)\left( {r\ln r + (r + 1)\ln \left( {\frac{2}{r + 1}} \right)}
\right)}
\right.\\
& \qquad \qquad \left. { - (1 - r)\left( {R\ln R + (R + 1)\ln \left(
{\frac{2}{R + 1}} \right)} \right)} \right].\notag
\end{align}
\end{example}

\begin{example} \label{ex310} (\textit{arithmetic-geometric divergence}). Let us
consider
\begin{equation}
\label{eq120}
f_T (x) = \left( {\frac{x + 1}{2}} \right)\ln \left( {\frac{x + 1}{2\sqrt x
}} \right),
\quad
x \in (0,\infty ),
\end{equation}

\noindent in (\ref{eq33}), then we have $C_f (P\vert \vert Q) =
T(P\vert \vert Q).$

Moreover,
\begin{equation}
\label{eq121}
{f}'_T (x) = \frac{1}{4}\left[ {1 - x^{ - 1} + 2\ln \left( {\frac{x +
1}{2\sqrt x }} \right)} \right],
\end{equation}

\noindent and
\begin{equation}
\label{eq122}
{f}''_T (x) = \frac{1}{4}\left( {\frac{1 + x^2}{x^2 + x^3}} \right).
\end{equation}

In view of (\ref{eq121}), (\ref{eq122}), Theorems \ref{the21} and
\ref{the22}, we have the following bounds on
\textit{arithmetic-geometric divergence}:
\begin{equation}
\label{eq123} 0 \leqslant T(P\vert \vert Q) \leqslant E_{T } (P\vert
\vert Q) \leqslant A_{T } (r,R)
\end{equation}

\noindent and
\begin{equation}
\label{eq124} 0 \leqslant T(P\vert \vert Q) \leqslant B_{T } (r,R)
\leqslant A_{T } (r,R),
\end{equation}

\noindent where
\begin{align}
\label{eq125} E_T (P\vert \vert Q) & = \frac{1}{4}\chi ^2(Q\vert
\vert P) + \frac{1}{2}\sum\limits_{i = 1}^n {(p_i - q_i )\ln \left(
{\frac{p_i + q_i }{2\sqrt {p_i q_i } }} \right)},\notag\\
A_T (r,R) & = \frac{1}{16}(R - r)^2\left[ {(rR)^{ - 1} + L_{ - 1}^{
- 1} \left( {r + 1,R + 1} \right) - L_{ - 1}^{ - 1} (r,R)}
\right]\notag\\
\intertext{and} B_T (r,R) & = \frac{1}{2(R - r)}\left[ {(R - 1)(r +
1)\ln \left( {\frac{r + 1}{2\sqrt r }} \right)}
\right.\notag\\
& \qquad \qquad \left. { + (1 - r)(R + 1)\ln \left( {\frac{R +
1}{2\sqrt R }} \right)} \right].\notag
\end{align}
\end{example}

\section{Relative J-Divergence and Inequalities}

In this section we shall present bound on \textit{relative
J-divergence} in terms of symmetric measures (\ref{eq13}),
(\ref{eq16})-(\ref{eq18}).

\begin{proposition} \label{pro41} (\textit{Relative J-divergence} and\textit{ triangular discrimination}). We have the following bounds:
\begin{equation}
\label{eq128}
\frac{(r + 1)(r + 3)}{8}\Delta (P\vert \vert Q) \leqslant D(P\vert \vert Q)
\leqslant \frac{(R + 1)(R + 3)}{8}\Delta (P\vert \vert Q),
\end{equation}
\end{proposition}

\begin{proof} Let us consider
\begin{equation}
\label{eq131} g_{D \Delta } (x) = \frac{{f}''_{D } (x)}{{f}''_\Delta
(x)} = \frac{(x + 1)(x + 3)}{8}, \quad x \in (0,\infty ),
\end{equation}

\noindent where ${f}''_{D } (x)$ and ${f}''_\Delta (x)$ are as given
by (\ref{eq50}) and (\ref{eq98}) respectively.

From (\ref{eq131}), we have
\begin{equation}
\label{eq132} {g}'_{D \Delta } (x) = \frac{2 + x}{4} > 0, \quad x
\in (0,\infty ).
\end{equation}

In view of (\ref{eq132}), we conclude that
\begin{equation}
\label{eq133} m = \mathop {\inf }\limits_{x \in [r,R]} g_{D \Delta }
(x) = \frac{(r + 1)(r + 3)}{8}
\end{equation}

\noindent and
\begin{equation}
\label{eq134} M = \mathop {\sup }\limits_{x \in [r,R]} g_{D \Delta }
(x) = \frac{(R + 1)(R + 3)}{8}.
\end{equation}

Expressions (\ref{eq133}) and (\ref{eq134}) together with
(\ref{eq43}) give the required result.
\end{proof}

\begin{proposition} \label{pro43} (\textit{Relative J-divergence} and\textit{ J-divergence}). We have the following bounds:
\begin{equation}
\label{eq142} \frac{r^2(r + 3)}{(r + 1)^3}J(P\vert \vert Q)
\leqslant D(P\vert \vert Q) \leqslant \frac{R^2(R + 3)}{(R +
1)^3}J(P\vert \vert Q),
\end{equation}
\end{proposition}

\begin{proof} Let us consider
\begin{equation}
\label{eq145} g_{D J} (x) = \frac{{f}''_{D } (x)}{{f}''_J (x)} =
\frac{x^2(x + 3)}{(x + 1)^3}, \quad x \in (0,\infty ),
\end{equation}

\noindent where ${f}''_{D } (x)$ and ${f}''_J (x)$ are as given by
(\ref{eq50}) and (\ref{eq106}) respectively.

From (\ref{eq145}), we have
\begin{equation}
\label{eq146} {g}'_{D J} (x) = \frac{6x}{(x + 1)^4} > 0, \quad x \in
(0,\infty ).
\end{equation}

In view of (\ref{eq146}), we conclude that
\begin{equation}
\label{eq147} m = \mathop {\inf }\limits_{x \in [r,R]} g_{D J} (x) =
\frac{r^2(r + 3)}{(r + 1)^3}
\end{equation}

\noindent and
\begin{equation}
\label{eq148} M = \mathop {\sup }\limits_{x \in [r,R]} g_{D J} (x) =
\frac{r^2(r + 3)}{(r + 1)^3}.
\end{equation}

Expressions (\ref{eq147}) and (\ref{eq148}) together with
(\ref{eq43}) give the required result.
\end{proof}

\begin{proposition} \label{pro45} (\textit{Relative J-divergence} and\textit{ Jensen-Shannon divergence}). We have the following bounds:
\begin{equation}
\label{eq156}
\frac{2r(r + 3)}{r + 1}I(P\vert \vert Q) \leqslant D(P\vert \vert Q)
\leqslant \frac{2R(R + 3)}{R + 1}I(P\vert \vert Q),
\end{equation}
\end{proposition}

\begin{proof} Let us consider
\begin{equation}
\label{eq159} g_{D I} (x) = \frac{{f}''_{D } (x)}{{f}''_I (x)} =
\frac{2x(x + 3)}{x + 1}, \quad x \in (0,\infty ),
\end{equation}

\noindent where ${f}''_{D } (x)$ and ${f}''_I (x)$ are as given by
(\ref{eq50}) and (\ref{eq114}) respectively.

From (\ref{eq159}), we have
\begin{equation}
\label{eq160} {g}'_{D I} (x) = \frac{2(x^2 + 2x + 3)}{(x + 1)^2}
> 0, \quad x \in (0,\infty ).
\end{equation}

In view of (\ref{eq160}), we conclude that
\begin{equation}
\label{eq161} m = \mathop {\inf }\limits_{x \in [r,R]} g_{D I} (x) =
\frac{2r(r + 3)}{r + 1}
\end{equation}

\noindent and
\begin{equation}
\label{eq162} M = \mathop {\sup }\limits_{x \in [r,R]} g_{D I} (x) =
\frac{2R(R + 3)}{R + 1}.
\end{equation}

Expressions (\ref{eq161}) and (\ref{eq162}) together with
(\ref{eq43}) give the required result.
\end{proof}

\section{Relative Jensen-Shannon divergence and Inequalities}

In this section we shall present bound on \textit{relative
Jensen-Shannon divergence} in terms of symmetric measures
(\ref{eq13}), (\ref{eq16})-(\ref{eq18}).

\begin{proposition} \label{pro51} (\textit{Relative Jensen-Shannon divergence}
and\textit{ triangular discrimination}). We have the following bounds:
\begin{equation}
\label{eq170}
\frac{R + 1}{8R}\Delta (P\vert \vert Q) \leqslant F(P\vert \vert Q)
\leqslant \frac{r + 1}{8r}\Delta (P\vert \vert Q),
\end{equation}
\end{proposition}

\begin{proof} Let us consider
\begin{equation}
\label{eq173} g_{F \Delta } (x) = \frac{{f}''_{F } (x)}{{f}''_\Delta
(x)} = \frac{x + 1}{8x}, \, x \in (0,\infty ),
\end{equation}

\noindent where ${f}''_{F } (x)$ and ${f}''_\Delta (x)$ are as given
by (\ref{eq66}) and (\ref{eq98}) respectively.

From (\ref{eq173}), we have
\begin{equation}
\label{eq174} {g}'_{F \Delta } (x) = - \frac{1}{8x^2} < 0, \quad x
\in (0,\infty ).
\end{equation}

In view of (\ref{eq174}), we conclude that
\begin{equation}
\label{eq175} m = \mathop {\inf }\limits_{x \in [r,R]} g_{F \Delta }
(x)  = \frac{R + 1}{8R}
\end{equation}

\noindent and
\begin{equation}
\label{eq176} M = \mathop {\sup }\limits_{x \in [r,R]} g_{F \Delta }
(x)  = \frac{r + 1}{8r}.
\end{equation}

Expressions (\ref{eq175}) and (\ref{eq176}) together with
(\ref{eq43}) give the required result.
\end{proof}

\begin{proposition} \label{pro53} (\textit{Relative Jensen-Shannon divergence} and\textit{ J-divergence}). We have the following bounds:
\begin{equation}
\label{eq184} 0 \leqslant F(P\vert \vert Q) \leqslant
\frac{4}{27}J(P\vert \vert Q),
\end{equation}
\end{proposition}

\begin{proof} Let us consider
\begin{equation}
\label{eq187} g_{F J} (x) = \frac{{f}''_{F } (x)}{{f}''_J (x)} =
\frac{x}{(x + 1)^3}, \quad x \in (0,\infty ),
\end{equation}

\noindent where ${f}''_{F } (x)$ and ${f}''_J (x)$ are as given by
(\ref{eq66}) and (\ref{eq106}) respectively.

From (\ref{eq187}), we have
\begin{equation}
\label{eq188} {g}'_{F J} (x) = - \frac{2x - 1}{(x +
1)^4}\begin{cases}
 { \geqslant 0,} & {x \leqslant \frac{1}{2}} \\
 { \leqslant 0,} & {x \geqslant \frac{1}{2}} \\
\end{cases},
\end{equation}

In view of (\ref{eq188}), we conclude that the function $g_{F J}
(x)$ is increasing in $(0,\frac{1}{2})$ and decreasing in
$(\frac{1}{2},\infty )$, and hence

\begin{equation}
\label{eq189} M = \mathop {\sup }\limits_{x \in [r,R]} g_{F J} (x) =
g_{F J} (\frac{1}{2}) = \frac{4}{27}.
\end{equation}

Now (\ref{eq189}) together with (\ref{eq43}) give the required
result.
\end{proof}

\begin{proposition} \label{pro55} (\textit{Adjoint of relative Jensen-Shannon divergence} and\textit{ Jensen-Shannon divergence}). We have the following bounds:
\begin{equation}
\label{eq198}
\frac{2}{R + 1}I(P\vert \vert Q) \leqslant F(P\vert \vert Q) \leqslant
\frac{2}{r + 1}I(P\vert \vert Q),
\end{equation}
\end{proposition}

\begin{proof} Let us consider
\begin{equation}
\label{eq201} g_{F I} (x) = \frac{{f}''_{F } (x)}{{f}''_I (x)} =
\frac{2}{x + 1}, \quad x \in (0,\infty ),
\end{equation}

\noindent where ${f}''_{F } (x)$ and ${f}''_I (x)$ are as given by
(\ref{eq66}) and (\ref{eq114}) respectively.

From (\ref{eq201}), we have
\begin{equation}
\label{eq202} {g}'_{F I} (x) = - \frac{2}{(x + 1)^2} < 0, \quad x
\in (0,\infty ).
\end{equation}

In view of (\ref{eq202}), we conclude that
\begin{equation}
\label{eq203} m = \mathop {\inf }\limits_{x \in [r,R]} g_{F I} (x) =
\frac{2}{R + 1}
\end{equation}

\noindent and
\begin{equation}
\label{eq204} M = \mathop {\sup }\limits_{x \in [r,R]} g_{F I} (x) =
\frac{2}{r + 1}.
\end{equation}

Expressions (\ref{eq203}) and (\ref{eq204}) together with
(\ref{eq43}) give the required result.
\end{proof}

\begin{remark} The inequalities (\ref{eq170})
and (\ref{eq198}) can also be as
\[
r \leqslant \zeta _t (P\vert \vert Q) \leqslant R, \quad t = 1 \,\,
and \,\,2,
\]

\noindent where
\[
\zeta _1 (P\vert \vert Q) = \frac{\Delta (P\vert \vert
Q)}{8F(P\vert \vert Q) - \Delta (P\vert \vert Q)},
\]

\noindent and
\[
\zeta _3 (P\vert \vert Q) = \frac{2I(P\vert \vert Q) - F(P\vert
\vert Q)}{F(P\vert \vert Q)}
\]

\noindent respectively.
\end{remark}

\section{Relative arithmetic-geometric divergence and Inequalities}

In this section we shall present bound on \textit{relative
arithmetic-geometric divergence} in terms of symmetric measures
(\ref{eq13}), (\ref{eq16})-(\ref{eq18}).

\begin{proposition} \label{pro61} (\textit{Relative arithmetic-geometric divergence }and\textit{ triangular discrimination}). We have the
following bounds:
\begin{equation}
\label{eq212}
\frac{(R + 1)^2}{16R^2}\Delta (P\vert \vert Q) \leqslant G(P\vert \vert Q)
\leqslant \frac{(r + 1)^2}{16r^2}\Delta (P\vert \vert Q),
\end{equation}
\end{proposition}

\begin{proof} Let us consider
\begin{equation}
\label{eq215} g_{G \Delta } (x) = \frac{{f}''_{G } (x)}{{f}''_\Delta
(x)} = \frac{(x + 1)^2}{16x^2}, \quad x \in (0,\infty ),
\end{equation}

\noindent where ${f}''_{G } (x)$ and ${f}''_\Delta (x)$ are as given
by (\ref{eq66}) and (\ref{eq98}) respectively.

From (\ref{eq215}), we have
\begin{equation}
\label{eq216} {g}'_{G \Delta } (x) = - \frac{x + 1}{8x^3} < 0, \quad
x \in (0,\infty ).
\end{equation}

In view of (\ref{eq216}), we conclude that
\begin{equation}
\label{eq217} m = \mathop {\inf }\limits_{x \in [r,R]} g_{G \Delta }
(x)  = \frac{(R + 1)^2}{16R^2}
\end{equation}

\noindent and
\begin{equation}
\label{eq218} M = \mathop {\sup }\limits_{x \in [r,R]} g_{G \Delta }
(x)  = \frac{(r + 1)^2}{16r^2}.
\end{equation}

Expressions (\ref{eq217}) and (\ref{eq218}) together with
(\ref{eq43}) give the required result.
\end{proof}

\begin{proposition} \label{pro63} (\textit{Relative arithmetic-geometric divergence }and\textit{ J-divergence}). We have the following
bounds:
\begin{equation}
\label{eq226} \frac{1}{2(R + 1)^2}J(P\vert \vert Q) \leqslant
G(P\vert \vert Q) \leqslant \frac{1}{2(r + 1)^2}J(P\vert \vert Q),
\end{equation}
\end{proposition}

\begin{proof} Let us consider
\begin{equation}
\label{eq229} g_{G J} (x) = \frac{{f}''_{G } (x)}{{f}''_J (x)} =
\frac{1}{2(x+1)^2}, \quad x \in (0,\infty ),
\end{equation}

\noindent where ${f}''_{G } (x)$ and ${f}''_J (x)$ are as given by
(\ref{eq66}) and (\ref{eq106}) respectively.

From (\ref{eq229}), we have
\begin{equation}
\label{eq230} {g}'_{G J} (x) = - \frac{1}{(x+1)^3} < 0, \quad x \in
(0,\infty ).
\end{equation}

In view of (\ref{eq216}), we conclude that
\begin{equation}
\label{eq231} m = \mathop {\inf }\limits_{x \in [r,R]} g_{G J} (x) =
\frac{1}{2(R + 1)^2}
\end{equation}

\noindent and
\begin{equation}
\label{eq232} M = \mathop {\sup }\limits_{x \in [r,R]} g_{G J} (x) =
\frac{1}{2(r + 1)^2}.
\end{equation}

Expressions (\ref{eq231}) and (\ref{eq232}) together with
(\ref{eq43}) give the required result.
\end{proof}

\begin{proposition} \label{pro65} (\textit{Relative arithmetic-geometric divergence }and\textit{ Jensen-Shannon divergence}). We have the
following bounds:
\begin{equation}
\label{eq240}
\frac{1}{R}I(P\vert \vert Q) \leqslant G(P\vert \vert Q) \leqslant
\frac{1}{r}I(P\vert \vert Q),
\end{equation}
\end{proposition}

\begin{proof} Let us consider
\begin{equation}
\label{eq243} g_{G I} (x) = \frac{{f}''_{G } (x)}{{f}''_I (x)} =
\frac{1}{x}, \quad x \in (0,\infty ),
\end{equation}

\noindent where ${f}''_{G } (x)$ and ${f}''_I (x)$ are as given by
(\ref{eq66}) and (\ref{eq114}) respectively.

From (\ref{eq243}), we have
\begin{equation}
\label{eq244} {g}'_{G I} (x) = - \frac{1}{x^2} < 0, \quad x \in
(0,\infty ).
\end{equation}

In view of (\ref{eq244}), we conclude that
\begin{equation}
\label{eq245} m = \mathop {\inf }\limits_{x \in [r,R]} g_{G I} (x) =
\frac{1}{R}
\end{equation}

\noindent and
\begin{equation}
\label{eq246} M = \mathop {\sup }\limits_{x \in [r,R]} g_{G I} (x) =
\frac{1}{r}.
\end{equation}

Expressions (\ref{eq245}) and (\ref{eq246}) together with
(\ref{eq43}) give the required result.
\end{proof}

\begin{proposition} \label{pro67} (\textit{Relative arithmetic-geometric divergence }and \textit{ arithmetic-geometric divergence}). We have the
following bounds:
\begin{equation}
\label{eq254} \frac{2}{1+R^2}T(P\vert \vert Q) \leqslant G(P\vert
\vert Q) \leqslant \frac{2}{1+r^2}T(P\vert \vert Q),
\end{equation}
\end{proposition}

\begin{proof} Let us consider
\begin{equation}
\label{eq257} g_{G T} (x) = \frac{{f}''_{G } (x)}{{f}''_T (x)} =
\frac{2}{1+x^2}, \quad x \in (0,\infty ),
\end{equation}

\noindent where ${f}''_{G } (x)$ and ${f}''_T (x)$ are as given by
(\ref{eq66}) and (\ref{eq122}) respectively.

From (\ref{eq257}), we have
\begin{equation}
\label{eq258} {g}'_{G T} (x) = - \frac{4x}{(1+x^2)^2} < 0, \quad x
\in (0,\infty ).
\end{equation}

In view of (\ref{eq258}), we conclude that
\begin{equation}
\label{eq259} m = \mathop {\inf }\limits_{x \in [r,R]} g_{G T} (x) =
\frac{2}{1+R^2}
\end{equation}

\noindent and
\begin{equation}
\label{eq260} M = \mathop {\sup }\limits_{x \in [r,R]} g_{G T} (x) =
\frac{2}{1+r^2}.
\end{equation}

Expressions (\ref{eq259}) and (\ref{eq260}) together with
(\ref{eq43}) give the required result.
\end{proof}

\begin{remark} The inequalities (\ref{eq212}),
(\ref{eq226}), (\ref{eq240}) and (\ref{eq254}) can also be as
\[
r \leqslant \xi _t (P\vert \vert Q) \leqslant R, \quad t = 1,2,3
\,\, and \,\,4,
\]

\noindent where
\[
\xi _1 (P\vert \vert Q) = \frac{\sqrt {\Delta (P\vert \vert Q)}
}{4\sqrt {G(P\vert \vert Q)} - \sqrt {\Delta (P\vert \vert Q)} },
\]
\[
\xi _2 (P\vert \vert Q) = \frac{\sqrt {J(P\vert \vert Q)} - \sqrt
{2G(P\vert \vert Q)} }{\sqrt {2G(P\vert \vert Q)} },
\]
\[
\xi _2 (P\vert \vert Q) = \frac{I(P\vert \vert Q)}{G(P\vert \vert
Q)},
\]

\noindent and
\[
\xi _4 (P\vert \vert Q) = \frac{\sqrt {2T(P\vert \vert Q) - G(P\vert
\vert Q)} }{\sqrt {G(P\vert \vert Q)} }
\]

\noindent respectively.
\end{remark}


\begin{thebibliography}{99}
\bibitem{pea} K. PEARSON, On
the Criterion that a given system of deviations from the probable
in the case of correlated system of variables is such that it can
be reasonable supposed to have arisen from random sampling,
\textit{Phil. Mag.}, \textbf{50}(1900), 157-172.
\end{thebibliography}
\end{document}